\documentclass[12pt]{amsart}
\usepackage{pstricks,pst-plot}

\parskip=\smallskipamount

\newtheorem{theorem}{Theorem}[section]

\newtheorem{proposition}[theorem]{Proposition}

\theoremstyle{definition}

\newtheorem{remark}[theorem]{Remark}

\numberwithin{equation}{section}

\newcommand{\B}{\mathbb{B}}
\newcommand{\C}{\mathbb{C}}

\newcommand{\Z}{\mathbb{Z}}
\renewcommand{\P}{\mathbb{P}}
\newcommand{\R}{\mathbb{R}}

\newcommand{\I}{\mathfrak{i}}

\newcommand{\cC}{\mathcal{C}}

\newcommand{\cL}{\mathcal{L}}

\newcommand\wt{\widetilde}
\newcommand\hra{\hookrightarrow}

\newcommand\di{\partial}

\newcommand\dist{\mathrm{dist}}
\newcommand\dibar{\overline\partial}

\newcommand{\ol}{\overline}

\usepackage{amssymb}
\input xy
\xyoption{all}

\begin{document}
\title[A complex surface]{A complex surface admitting a strongly plurisubharmonic function but no holomorphic functions}
\author{Franc Forstneri\v c}
\address{F.\ Forstneri\v c, Faculty of Mathematics and Physics, University of Ljubljana, and Institute of Mathematics, Physics and Mechanics, Jadranska 19, 1000 Ljubljana, Slovenia}
\email{franc.forstneric\,@fmf.uni-lj.si}

%
%
\date{October 29, 2012}

\begin{abstract}
We find a domain $X\subset \C\P^2$ with a strongly plurisubharmonic function such that every holomorphic function on $X$ is constant. 
\end{abstract}

\maketitle

\section{The main result}
\label{sec:Intro}
One of the most important classes of functions in complex analysis is the class of {\em strongly plurisubharmonic} functions.
By a classical result of Grauert \cite{Grauert:Levi}, a complex manifold with an exhausting strongly plurisubharmonic function is a {\em Stein manifold}, i.e., it admits plenty of holomorphic functions. It is therefore natural to ask whether the existence of a (non-exhausting) strongly plurisubharmonic function on a complex manifold $X$ tells us anything at all about holomorphic functions on $X$. In particular, does such $X$ necessarily admit a noncostant holomorphic function? Is $X$ a domain in a Stein manifold? 

This apparently long standing problem was communicated to me by Karl Oeljeklaus on Oct.\ 25, 2012 at a conference in Lille. (According to F.\ L\'arusson, the question was known to R.\ Narasimhan more than 20 years ago.) The purpose of this note is to explain a class of counter-examples which are domains in the projective plane $\C\P^2$. Another example, also in $\C\P^2$, was proposed by T.-C.\ Dinh. One can find similar examples in any $\C\P^n$ for $n>1$; see Remark \ref{rem:higherdim} below due to S.\ Nemirovski.

\begin{theorem}
\label{th:main}
There exists a connected open set $X$ in $\C\P^2$ such that $X$ admits a strongly plurisubharmonic function, but every holomorphic function on $X$ is constant.
\end{theorem}

We can choose $X$ to be a small tubular neighborhood of a generically embedded closed oriented real surface $S\hra \C\P^2$ of genus $g(S)\le 2$ and of degree $d\ne 0$. In particular, a generic perturbation of the projective line $\C\P^1 \subset \C\P^2$ admits such neighborhoods.

Our example is a consequence of the following two results. The first of them is due to S.\ Nemirovski \cite{Nemirovski1999}.

\begin{theorem}
\label{th:Nem}
If $S$ is closed oriented real surface of genus $g\le 2$, smoothly embedded in the projective plane $\C\P^2$ with $0\ne [S]\in H_2(\C\P^2;\Z)$, then every holomorphic function in an open connected neighborhood of $S$ is constant. This holds in particular for any smoothly embedded homologically nontrivial 2-sphere in $\C\P^2$. 
\end{theorem}

On the other hand, every closed oriented surface $S$ of genus $g\ge 3$ admits a smooth embedding $S\hra \C\P^2$ with a basis of open Stein neighborhoods in $\C\P^2$. In fact, every embedding $S\hra \C\P^2$ of degree $d>0$ satisfying $2g \ge d^2 + 3d + 2$ admits a smooth isotopy to an embedding with a Stein neighborhood basis; see Corollary 9.6.2 in \cite{F:book}.  

The main ingredient in the proof of Theorem \ref{th:Nem} is the {\em generalized adjunction inequality} (see (\ref{eq:GAI}) below) which depends on highly nontrivial results from Seiberg-Witten theory. This interesting connection was observed and explored independently by Nemirovski \cite{Nemirovski1999} and by Lisca and Mati\'c \cite{Lisca-Matic1997}.

Our second ingredient is completely elementary. Let $X$ be a complex manifold. Denote by $J$ the almost complex structure operator induced by the complex structure on $X$; thus $J$ is a real endomorphism of the tangent bundle $TX$ satisfying $J^2=-Id$. A real submanifold $S$ of class $\cC^1$ in $X$ is said to be {\em totally real} at a point $p\in S$ if the tangent space $T_p S$ satisfies $T_p S\cap J(T_p S)=\{0\}$; $S$ is totally real if this holds at each point. A point $p\in S$ which is not totally real is said to be a {\em complex point} of $S$.

\begin{proposition}
\label{prop:SPSH}
Let $S$ be a real submanifold of class $\cC^1$ embedded in a complex manifold $X$. If $S$ has only isolated complex points, then there exists a strongly plurisubharmonic function in an open neighborhood of $S$ in $X$.
\end{proposition}

Proposition \ref{prop:SPSH} is proved in Sect.\ \ref{sec:SPSH} below, and is a special case of Proposition \ref{prop:SPSH2}. We hope that this may be of some independent interest. 

\smallskip
\noindent{\it Proof of Theorem \ref{th:main}.} 
Let $S$ be an oriented closed real 2-surface, smoothly embedded in $\C\P^2$, with nontrivial homology class $0\ne [S] \in H_2(\C\P^2;\Z)$. After a generic deformation of $S$ we may assume that $S$ only contains isolated complex points. (For instance, there is a small perturbation of the projective line $\C\P^1\subset \C\P^2$ with 3 elliptic complex points and no hyperbolic complex points \cite{FF:complexpoints}, but this additional information will not be used.) By Proposition \ref{prop:SPSH} there exists a connected open neighborhood $X\subset \C\P^2$ of $S$ and a strongly plurisubharmonic function $\rho\colon X\to\R$. On the other hand, if the genus of $S$ satisfies $g<3$, then Theorem \ref{th:Nem} shows that $X$ does not admit any nonconstant holomorphic functions. 
\qed
\smallskip

For completeness we include a sketch of proof of Theorem \ref{th:Nem}; for details we refer to \cite{Nemirovski1999}, or to the survey in \cite[Chap.\ 9]{F:book}; see in particular Corollary 9.8.2 in \cite{F:book}. Let $S\hra \C\P^2$ be as in Theorem \ref{th:Nem}. Recall that $H_2(\C\P^2;\Z)\cong\Z$, the generator being the class of the projective line. Using this identification, let $[S]=d\ne 0$, so $d\in\Z$ is the degree of $S$. Suppose that there exists a nonconstant holomorphic function in an open connected neighborhood $U$ of $S$. The envelope of holomorphy of $U$ is then a Stein Riemann domain $\Omega\to\C\P^2$  which contains an embedded copy of a neighborhood of $S$ in $\C\P^2$. (See the discussion and references in \cite[\S 1.4]{Nemirovski1999}.) Clearly $S$ remains homologically nontrivial in $\Omega$. Any relatively compact subdomain of $\Omega$ embeds as a domain in a compact K\"ahler surface $Y$ with $b_2^+(Y)>1$ and with ample canonical bundle; this embeds $S$ to $Y$. Now the results of Seiberg-Witten theory imply the following {\em generalized adjunction inequality}:
\begin{equation}
\label{eq:GAI}
	2g(S)-2 \ge [S]\cdotp [S]  + \bigl| K_Y\cdotp [S] \bigr|.
\end{equation}
(See \cite[Theorem 9.7.1]{F:book} for a summary of these results.) Here $[S]\cdotp[S]$ is the self-intersection number of $S$, and $K_Y \cdotp [S]$ is the value on $[S]$ of the canonical class $K_Y$ of $Y$ (the negative $-c_1(TY)$ of the first Chern class of the tangent bundle $TY$). Both numbers on the right hand side of (\ref{eq:GAI}) can be computed by replacing $Y$ with any domain containing $S$; hence they equal the corresponding numbers with $Y$ replaced by $\C\P^2$. This gives $[S]\cdotp [S]=d^2$, $|K_Y\cdotp [S]|=3|d|$, and hence $2g(S)\ge d^2 + 3|d| + 2\ge 6$, so $g(S)\ge 3$.

\begin{remark}
The simplest complex surface $X$ satisfying Theorem \ref{th:main} that we obtain by the above argument is diffeomorphic to the normal bundle of the projective line in $\C\P^2$. The zero section of this bundle is an embedded sphere with self-intersection number $+1$. A smooth oriented 4-manifold containing such a sphere does not admit any Stein structure since a homologically nontrivial embedded 2-sphere $S$ in any Stein surface satisfies $[S]\cdotp [S]\le -2$ as follows from (\ref{eq:GAI}). On the other hand, since our $X$ is a 2-dimensional CW-complex, it admits a non-tame Stein structure which is realizable as a domain in $\C^2$ (see Gompf \cite{Gompf1,Gompf2} and \cite[Sec.\ 9.12]{F:book}).
\end{remark}

\begin{remark}
\label{rem:higherdim} 
After disseminating an earlier version of this paper, Stefan Nemirovski observed that the same argument also gives such examples in higher dimensions (personal communication). Let us begin with a homologically nontrivial embedded 2-sphere $S\subset \C\P^2$ with isolated complex points. Now embed $\C\P^2$ holomorphically (for instance, linearly) into $\C\P^n$ for some $n>2$. Proposition \ref{prop:SPSH} shows that $S$ has a connected open neighborhood $U\subset \C\P^n$ with a strongly plurisubharmonic function on it. Every holomorphic function on $U$ is constant on $U\cap \C\P^2$ by the above argument. The envelope of holomorphy $\wt U\to\C\P^n$ of $U$ is either a Stein domain over $\C\P^n$, or $\C\P^n$ itself (see the discussion and references in \cite[\S 1.4]{Nemirovski1999}). Since $U$ is schlicht, it embeds into its envelope $\wt U$. Hence, since holomorphic functions do not separate points in $U$, the envelope $\wt U$ equals $\C\P^n$ by the above alternative, so all holomorphic functions on $U$ are constant.

Note that this example is highly unstable --- we can deform the embedded sphere $S\subset \C\P^n$ for $n>2$ to a totally real embedded sphere which then admits a basis of open Stein neighborhoods.
\end{remark}


\section{Strongly plurisubharmonic functions near stratified totally real sets} 
\label{sec:SPSH}
In this section we prove Proposition \ref{prop:SPSH} and indicate a generalization to the class of {\em stratified totally real sets}; see Proposition \ref{prop:SPSH2} below. 

We begin with some preliminaries. Recall that the {\em Levi form}, $\cL_\rho$, of a $\cC^2$ function $\rho$ on a complex manifiold $X$ is a quadratic Hermitean form on the tangent bundle $TX$ that is given in any local holomorphic coordinate system $z=(z_1,\ldots,z_n)$ on $X$ by the {\em complex Hessian}:
\[
   \cL_\rho(z;v)=\sum_{j,k=1}^n \frac{\di^2 \rho(z)}{\di z_j\di\bar z_k}\, v_j\ol v_k,\quad v\in \C^n.
\]
(We identify the real tangent bundle $TX$ in the standard way with the holomorphic tangent bundle $T^{1,0} X$, the $(1,0)$-part of $TX\otimes \C$; see \cite[Sec.\ 1.6]{F:book}. For an intrinsic definition of the Levi form, associated to the $(1,1)$-form $dd^c\rho=\I \di\dibar\rho$, see e.g.\ \cite[Sec.\ 1.8]{F:book}.) The function $\rho$ is plurisubharmonic if $\cL_\rho\ge 0$, and is strongly plurisubharmonic if $\cL_\rho>0$; clearly the latter is an open condition.

The following well known results can be found in the paper \cite{Range-Siu1974} by Range and Siu; an exposition is available in several other sources, for instance in \cite[Theorem 6.1.6]{Stout:polynomial} and in \cite[Sec.\ 3.5]{F:book}. 

Let $M$ be a locally closed totally real submanifold of a complex manifold $X$. If $\tau \ge 0$ is a nonnegative $\cC^2$ function in an open neighborhood of $M$ which vanishes exactly to order two on $M$ (in the sense that its real Hessian is positive definite in the normal direction to $M$), then $\tau$ is strongly plurisubharmonic along $M$. If $M$ is of class $\cC^2$ then the function $\tau(x)=\dist(x,M)^2$, the squared distance to $M$ with respect to any smooth Riemannian metric on $X$, is such; if $M$ is merely of class $\cC^1$, we can apply Whitney's jet extension theorem to find a function $\tau$ with these properties. From this and the theorem of Grauert \cite{Grauert:Levi} it follows easily that every totally real submanifold $M$ admits a basis of open Stein neighborhoods in $X$; see \cite{Range-Siu1974} or \cite[Corollary 3.5.2]{F:book}. 

We make another trivial but useful observation: If $\rho$ vanishes to the second order along a submanifold $M\subset X$ and $h$ is a $\cC^2$ function in a neighborhood of $M$, then we have
\[
	\cL_{h\rho}(x;v)= h(x) \cL_\rho(x;v),\quad  x\in M,\ v\in T_x X. 
\]
Indeed, the Levi form $\cL_{h\rho}$ contains second order partial derivatives of the product $h\rho$; since $\rho$ vanishes to second order on $M$, the conclusion follows from the Leibniz formula.

\smallskip
\noindent{\em Proof of Proposition \ref{prop:SPSH}.} 
Let $S\subset X$ be a $\cC^1$ submanifold which is totally real, except at a discrete set of points $P=\{p_j\}\subset S$. By what has been said above, there exists a strongly plurisubharmonic function $\tau\ge 0$ in a neighborhood of $M:=S\setminus P$ in $X$ that vanishes to second order on $M$.

For every index $j$ we choose an open neighborhood $U^j\subset X$ of the point $p_j$ and a holomorphic coordinate map 
\[
	\theta^j=(\theta^j_1,\ldots,\theta^j_n)\colon U^j \stackrel{\cong}{\longrightarrow} \B \subset \C^n
\] 
onto the unit ball $\B=\{z\in\C^n \colon |z|<1\}$, with $n=\dim X$. We may assume that the sets $\overline U^j$ are pairwise disjoint. We denote by $r\B \subset \C^n$ the ball or radius $r>0$, and we set $U^j_r=(\theta^j)^{-1}(r\B) \subset U^j$ for any $0<r\le 1$. 

Pick a smooth function $\xi \colon X\to [0,1]$ such that $\xi=1$ on $X\setminus \cup_j U^j_{1/3}$, and $\xi=0$ in a small neighborhood of any point $p_j$. The function $\xi\tau$ is then defined in a full open neighborhood of $S$ in $X$, and its Levi form satisfies 
\begin{equation*}
\label{eq:e1}
	\cL_{\xi \tau}(x,\cdotp)=\xi(x)\cL_{\tau}(x;\cdotp)\ge 0,\qquad x\in S.
\end{equation*}
At points $x\in S\setminus\cup_j U^j_{1/3}$ we have $\xi=1$ and hence $\cL_{\xi \tau}(x;\cdotp) = \cL_{\tau}(x;\cdotp)>0$.

Choose a smooth function $\chi\colon \R\to [0,1]$ such that $\chi(t)=1$ for $t\le 1/2$ and $\chi(t)=0$ for $t\ge 2/3$. Set $|\theta^j|^2=\sum_{k=1}^n |\theta^j_k|^2$. The smooth function 
\[
	\phi_j:=\chi(|\theta^j|) \cdotp |\theta^j|^2 \colon U_j\longrightarrow \R_+
\]
agrees with $|\theta^j|^2$ on $U^j_{1/2}$, so it is strongly plurisubharmonic there, and its support is contained in $\ol U^j_{2/3}$. We extend $\phi_j$ to $X$ as the zero function on $X\setminus U^j_{2/3}$. Pick a constant $c_j>0$ such that the Levi form of the function $\xi\tau+c_j\phi_j$ is positive definite at each point of $x\in S\cap U^j$. In fact, this holds for every sufficiently small $c_j>0$
since $\cL_{\xi\tau}(x;\cdotp)\ge 0$ for every such $x$, and $\cL_{\xi\tau}(x;\cdotp) = \cL_{\tau}(x;\cdotp)> 0$ when $x\in U^j\setminus U^j_{1/2}$. It follows that the function
\[
	\rho:= \xi \tau+ \sum_j c_j \phi_j 
\]
is defined in a neighborhood of $S$ in $X$ (the sum is locally finite since the functions $\phi_j$ have pairwise disjoint supports), and it satisfies $\cL_\rho(x;\cdotp)>0$ for every point $x\in S$. Thus $\rho$ is strongly plurisubharmonic in a neighborhood of $S$ in $X$. This proves Proposition \ref{prop:SPSH}.
\qed
\smallskip

The above proof easily generalizes to the following situation.

\begin{proposition}
\label{prop:SPSH2}
Let $S$ be a locally closed subset of a complex manifold $X$ which admits a stratification 
$S=S_0\supset S_1\supset \cdots\supset S_m \supset S_{m+1}=\varnothing$ by locally closed sets $S_k$ such that every difference $S_k\setminus S_{k+1}$ is a union of pairwise disjoint $\cC^1$ totally real submanifolds of $X$. Then there exists a strongly plurisubharmonic function in an open neighborhood of $S$ in $X$.
\end{proposition}

A subset $S$ as in the above proposition could reasonably be called {\em stratified totally real set} in $X$.  It seems likey that a generic perturbation of any smooth submanifold $S\subset X$ of dimension $\dim_\R S\le \dim_\C X$ yields a stratified totally real submanifold, but we shall not pursue this issue here. 

To prove Proposition \ref{prop:SPSH2} we proceed as before by choosing a nonnegative strongly plurisubharmonic function in a neighborhood of the lowest stratum $S_m$ and vanishing on $S_m$; we then cut if off outside a suitable neighborhood of $S_m$ and add to it a strongly plurisubharmonic function in a neighborhood of the next stratum $S_{m-1}\setminus S_m$, etc. We leave the details to the reader.

\medskip\noindent
{\bf A final remark.} A reader familiar with Siu's theorem \cite{Siu1976} (that every Stein subvariety in a complex manifold admits a strongly plurisubharmonic function in its neighborhood, and also a Stein neighborhood basis) may notice some similarities between the two proofs, although the situation in this paper is substantially simpler to treat. An exposition of Siu's theorem, and of the extensions to the $q$-convex case given independently by Demailly and Col\c toiu, can be found in Sec.\ 3.1--3.2 of \cite{F:book}. 

Confronting these two results that hold under very different (essentially opposite) hypotheses, a natural question is what could be said in the mixed case. For example, assume that $S$ is a noncompact Levi-flat submanifold foliated by Stein complex leaves; does such $S$ admit a strongly plurisubharmonic function in a neighborhood of any compact subset $K\subset S$? The existence of a Stein neighborhood basis of Levi-flat submanifolds is a much more delicate issue as is seen on the example of the Diederich-Forn\ae ss {\rm worm domain} \cite{Diederich-Fornaess1977}. Perhaps one might expect positive results in the case of {\em simple} Levi-flat foliations of hypersurface type; see the related results in \cite{FL}.

\medskip
{\bf Acknowledgements.} 
I wish to acknowledge support by the research program P1-0291 from ARRS, Republic of Slovenia. I thank Karl Oeljeklaus for communicating the question answered in this note, and Stefan Nemirovski for providing the example in Remark \ref{rem:higherdim}. Last but not least, I thank the organizers of the conference {\em Geometric Methods in Several Complex Variables} (Lille, France, 25--26 October 2012) for their kind invitation and hospitality.

\bibliographystyle{amsplain}

\end{document}